\documentclass[10pt]{article}

\usepackage{graphicx}
\usepackage{algorithm}
\usepackage{algorithmic}
\usepackage{subfigure}

\usepackage[a4paper, left=35mm,right=35mm,top=34mm,bottom=34mm]{geometry}
\usepackage[utf8]{inputenc}
\usepackage[T1]{fontenc}
\usepackage[english]{babel}

\usepackage{fancyhdr}

\title{Point-Cloud-based Deep Learning Models for Finite Element Analysis}

\author{Meduri Venkata Shivaditya$^*$,
Francesca Bugiotti\thanks{Universit\'e Paris-Saclay, CNRS, CentraleSup\'elec, LISN, 1 rue Raimond Castaing, 91190 Gif-sur-Yvette, France},
Fr\'ed\'eric Magoul\`es\thanks{Universit\'e Paris-Saclay, CentraleSup\'elec, MICS, 9 rue Joliot Curie, 91192 Gif-sur-Yvette, France, Email: frederic.magoules@hotmail.com}
}

\begin{document}

\maketitle
\thispagestyle{fancy}

\begin{abstract}
In this paper, we explore point-cloud based deep learning models to analyze numerical simulations arising from finite element analysis.
The objective is to classify automatically the results of the simulations without tedious human intervention. 
Two models are here presented: the Point-Net classification model and the Dynamic Graph Convolutional Neural Net model.
Both trained point-cloud deep learning models performed well on experiments with finite element analysis arising from automotive industry.
The proposed models show promise in automatizing the analysis process of finite element simulations.
An accuracy of 79.17\% and 94.5\% is obtained for the Point-Net and the Dynamic Graph Convolutional Neural Net model respectively. 
\end{abstract}


\section{Introduction}

Finite element methods (FEM) are extensively used to simulate real world physical phenomenon by solving the associated partial differential equations \cite{Cia2002}, \cite{Hug2003}.
FEM have been applied to various fields including computational mechanics, fluid dynamics, and heat transfer for instance.
In FEM, the mathematical equations are reformulated with a variational formulation which is discretized in time and space.
Discretization is performed on a mesh, and the quality of the elements of the mesh impacts directly the approximated solution.
Several techniques exit to ensure high quality of the approximated solution, including for instance Streamline-Upwind Petrov-Galerkin (SUPG) \cite{Rus2006}, stabilized finite elements \cite{MZ2018CMAME}, bubble elements \cite{DF2008}.
The FEM solution is then obtained by solving a linear system of equations, whose size is proportional to the number of discretization points composing the mesh.
Simulations generated by FEM are then usually verified manually following visual inspection and preliminary checks.
This is a tedious process and requires human intervention.

The motivation behind this paper is to automate the process of analyzing the simulation results generated by FEM without human intervention. 
We explore the application of deep learning models to automatically analyze the simulations generated by FEM process.
This paper focuses on the task of classifying the generated FEM results based on a small number of human generated labels of so called 'normal' and 'abnormal' results.

In FEM, a mesh represents a domain into a discrete number of elements for which the solution can be calculated.
A mesh consists of geometrical information (nodes coordinates), a topology information (connectivity between nodes and elements), and field data stored on the mesh.
Opposite, Point-Clouds are simply a set of data points containing spatial coordinates and point features.
If only the nodes of the mesh and the fields on the mesh are considered, we can define a mapping between the mesh and a point cloud object.
Both mesh and Point-cloud objects share the property of permutation in-variance. 
Learning on point-cloud objects is researched by the scientific community and several deep learning models have been developed for several tasks including object classification and object segmentation \cite{WG2021}, \cite{SEG2021}.
Point-Net and Dynamic Graph Convolutional Neural Networks (DGCNN) \cite{SGT2008} \cite{YSL2019} are standard models for learning on point-cloud objects.

Finite element simulation of acoustic noise inside a car compartment has been used in this paper.
Depending on interior and exterior excitations sources, simulations show regions of high acoustic pressure near the front windows of the car compartment.
The results are divided into two classes: (i) class of results which have high noise level near the front windows, and (ii) class of results which have low noise level near the front windows.
The results are divided into classes by visual inspection by humans and a small dataset of 72 meshes with labels has been generated.  
The dataset is subdivided into training and testing sets on which the model is built and tested.
Point-Net and DGCNN models will be trained on the training meshes and labels.
After which, the trained models are evaluated using the accuracy metric on the test set. 

This paper shows that the state-of-the-art point-cloud classification models perform well for the task of classification of the FEM generated simulation.
The Point-Net model achieved an accuracy of 79.17\% while the DGCNN model was able to achieve an accuracy of 94.5\% on the test set.
As the models perform well, the point-cloud deep learning models can be used for automatic analysis of FEM generated simulation. 

This paper is organized as follows,
In section~\ref{section:relwor}, we present state-of-the-art models for the task of classification on point-cloud objects.
In section~\ref{section:method}, we present the methodology followed including framing the problem, exploring the dataset, evaluation methods.
Model architectures and results are presented in section~\ref{sec:result}.
In section~\ref{sec:conc}, we present the findings of the study. % and scope for future research on the subject.

\section{Related Works}
\label{section:relwor}

Classification is the problem of identifying which of a set of categories an object belongs to.
There are different types of algorithms to create a classification model, see for instance \cite{ZE2012RSER} and references therein.
The most popular of which are Support Vector Machines (SVM) \cite{HDO1998} \cite{ZE2012JACT}, Multi-layer perceptron, Convolutional Neural Networks, recursive deterministic perceptron (RDP) Neural Network, and Deep Learning Methods. 

SVMs construct a hyperplane that separates the objects with the maximum possible margin.
Taking advantage of the Kernel Trick, we are able to project the data into higher dimensions which allows for non-linear margins \cite{ZE2010JACT}. 

Multi-layer perceptrons (MLP), Convolutional Neural Networks (CNN) \cite{MZE2013EB}, Deep Learning models all are made up of neurons.
MLPs are stacked fully connected layers where a neuron of a layer is connected to every neuron in the other layer.
In CNNs, instead of matrix multiplications, convolutional operation is employed.
Due to increase in compute power, deeper models with millions of parameters are able to be trained by the computer which resulted in an increase in performance of the model.  

Several works in the research community study deep learning on point sets.
Point-Net \cite{QCH2017} is a highly efficient and effective point-cloud learning model and has been used as baseline for several point-sets learning research projects.
Before point-net, researchers used to transform point-cloud data into ND voxel grids or collections of images which renders data unnecessarily voluminous \cite{VOX2021}.
Point-Cloud models have been tested on the industry standard dataset for point clouds known as the ModelNet40.
The Point-Net model achieved an accuracy of 89.2\% on the ModelNet40 dataset. 

However, by design Point-Net does not capture local structures induced by the metric space points live in, limiting its ability to recognize fine-grained patterns and its generalization to complex scenes.
The DGCNN model \cite{YSL2019} allows the model to learn local level features by taking into account, the local neighbourhood using the $K$-Nearest Neighbours method.
The novel approach of edge-convolution and taking into advantage the local level features, the DGCNN model is much more powerful than the Point-Net model.
The DGCNN model achieved an accuracy of 92.9\% on the ModelNet40 dataset.
The DGCNN model involves computation of $K$-Nearest Neighbours after every layer, which makes the model time exhaustive and slows model training significantly compared to Point-Net model. 

These models are very efficient in learning the patterns of point-cloud objects.  These models have been employed in this paper to create the binary classification model of point-clouds generated from FEM meshes.

\section{Methodology}
\label{section:method}

\emph{Problem statement:}
We are interested by a function $F$ that takes the input point-cloud object as input and returns class label of the point-cloud object, where $F$ is a Binary Classification function of the form $F : N \times (d + C) \to Y$.
A point cloud is represented as a set of $N$ points $\{P_i| i = 1, \ldots, N\}$, where each point $P_i$ is a vector of $(x, y, z)$-coordinates plus extra feature channels such as acoustic pressure, etc.
The point-cloud is represented as a matrix $N \times (d + C)$, where $d$ is the dimension of the spatial coordinates and $C$ denotes the additional feature channels.
Quantity $Y$ denotes the class of the point-cloud object.

\begin{figure}[h]
\centering
\includegraphics[width=0.405\linewidth]{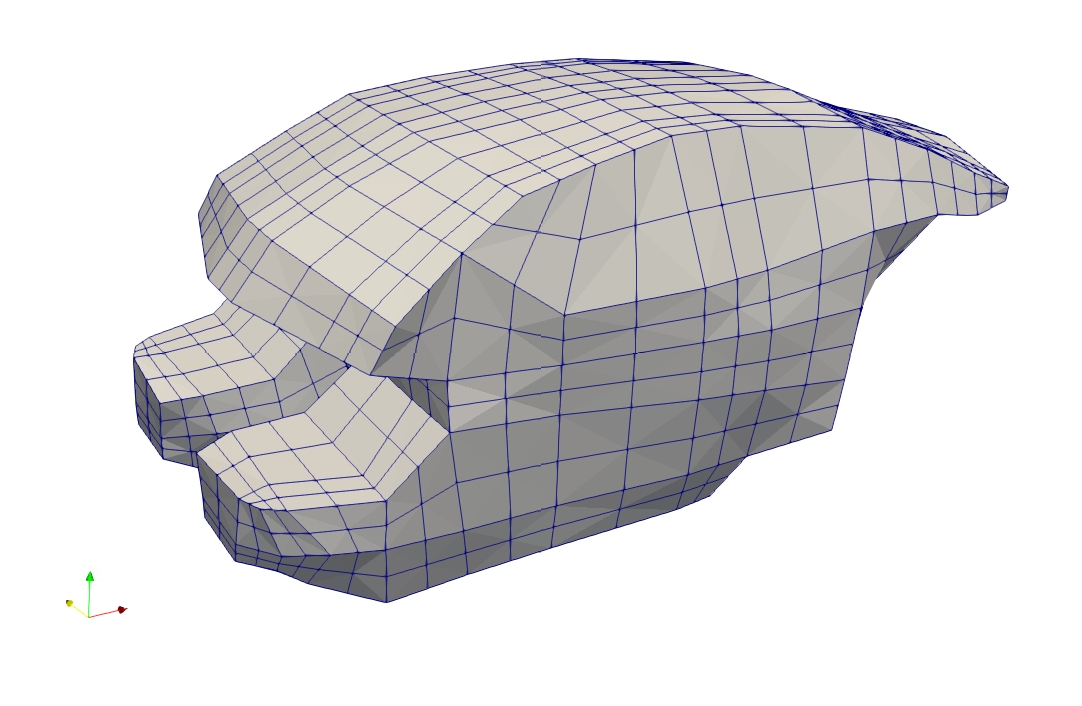}
\includegraphics[width=0.405\linewidth]{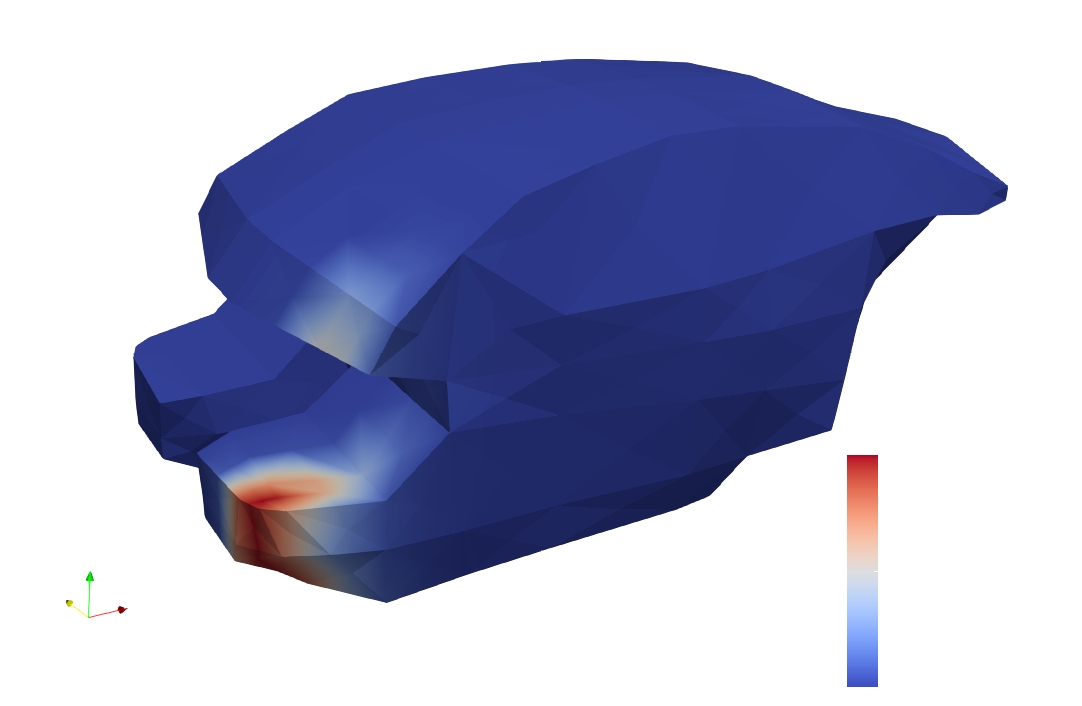}
\caption{Example of a finite element mesh (left) and acoustic pressure in Pa (right).}
\label{fig:solution}
\end{figure}

\emph{Dataset:}
The dataset is generated using FEM method.
An example of a finite element mesh is illustrated in Figure~\ref{fig:solution} (left) and the pressure obtained from the finite element method with Lagrange elements is represented in Figure~\ref{fig:solution} (right).
A total of 72 results have been generated and labels are given to each of them manually through visual inspection. 
The dataset is subdivided into training and testing sets as follow: 66.6\% of the dataset is used for training and 33.3\% of the dataset is used for testing.

\emph{Evaluation metric:}
Binary Classification in deep learning is evaluated generally using the Categorical Cross Entropy loss function.
In every iteration, the loss on the training set is computed and gradient descent is performed to modify the weights of the model in the direction of minimizing the Categorical Cross Entropy Loss, i.e.,
$$ -\sum_{c=1}^My_{o,c}\log(p_{o,c}) $$
where $M$ denotes the number of classes, in this case it is equal to 2, $y$ denotes binary indicator (0 or 1) if class label $c$ is the correct classification for observation $o$ and $p$ denotes the predicted probability observation $o$ is of class $c$.

\emph{Image classification models:}
The point-cloud data can be converted into ND voxel grids or collections of images.
The given data consists of 2D spatial coordinates and 1D additional pressure feature channel.
This data object can be converted into 2D gray-scale image following Algorithm~1. 

\begin{algorithm}
    \caption{Point Cloud to Projected Image on a plane}\label{trainalgo}
    \begin{algorithmic}
        \STATE  $\mathrm{pointcloudToImage} (PC(pointcloud), fields, plane)$
        \STATE  NDIM = DIM(plane)
        \STATE  NPOINTS = LENGTH(PC)
        \STATE  IMG = ZEROARRAY(NDIM)
        \FOR{point = $1$ to $NPOINTS$}
            \STATE IMG[PC[point][plane]] += fields[point]
        \ENDFOR
        \STATE return IMG
    \end{algorithmic}
\end{algorithm}

\begin{figure}[h]
    \centering
    \includegraphics[width=0.607\linewidth]{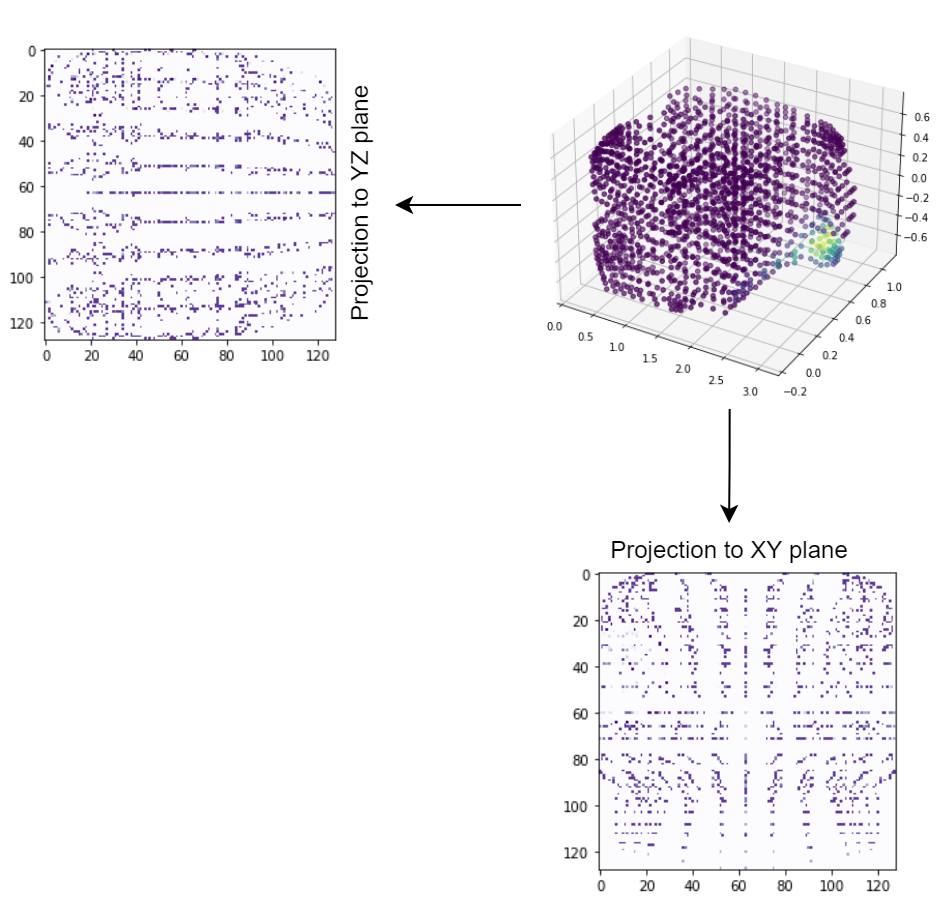}
    \caption{Image processing based on cutting plane.}
		\label{fig:image-processing}
\end{figure}

Now, the images can be directly used to train a Convolutional Neural Network (CNN).
A CNN model with three Convolutional Layers is defined for this model and the model is trained for 20-30 epochs on Categorical Cross Entropy Loss.
The model requires defining the data in ND voxel grid.
For 2D point-clouds, this is feasible in local memory.
For 3D point-clouds with high resolution, this model is infeasible to store in memory as the ND voxel image will require a lot of memory. In order to convert the 3D point-cloud to a 2D image, the 3D point-cloud should be projected on a plane which will be later considered for CNN model training.
The main drawback is to obtain the image, since it depends on the angle of view.
Figure~\ref{fig:image-processing} presents 3D point-cloud projection onto two planes, i.e., the $XY$ and the $YZ$ plane.

\emph{Point-cloud classification models:}
The models explored here are all belonging to the deep learning framework.
Deep Learning models are trained using Stochastic Gradient Descent where the weights are updated in the direction that minimizes the loss.
Algorithm~2 details the training algorithm considered, where $f$ represents the deep learning model, $y$ represents the actual class of the object.

\begin{algorithm}
    \label{algo:train}
    \caption{Training algorithm}\label{trainalgo}
    \begin{algorithmic}
        \STATE  $\mathrm{train\_model} (f, dataloader, EPOCHS)$
        \STATE Randomly initialize model weights $w$
        \FOR{epochs = $1$ to $EPOCHS$}
            \FOR{(pc, y) in dataloader} 
                \STATE input pc(Point-Cloud) object in the input layer
                \STATE forward propagate to get predicted result $\hat{y}$
                \STATE compute $e=CategoricalCrossEntropy(y, \hat{y})$
                \STATE back propagate $e$ from right to left through layers
                \STATE update $w$
            \ENDFOR
        \ENDFOR
    \end{algorithmic}
\end{algorithm}

\begin{figure}[h]
    \centering
    \includegraphics[width=0.95\linewidth]{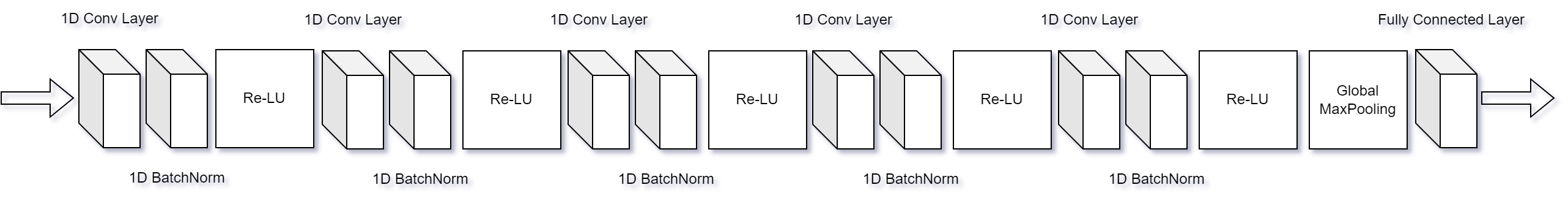}
    \caption{Point-net model architecture.}
		\label{fig:pointnet}
\end{figure}

\emph{Point-Net model:}
Point-Net model comprises of 1D convolution layers which project the point features into higher dimensions through the layers to learn features of the point-set.
The model also comprises of 1D batch normalization layers as well as Re-LU activation layers.
The last layer of the model is modified for the task of binary classification, i.e., the last layer of the model is a fully connected layer which outputs a vector of size 2.
The point-net model architecture is represented in Figure~\ref{fig:pointnet}.

\begin{figure}[h]
    \centering
    \includegraphics[width=0.95\linewidth]{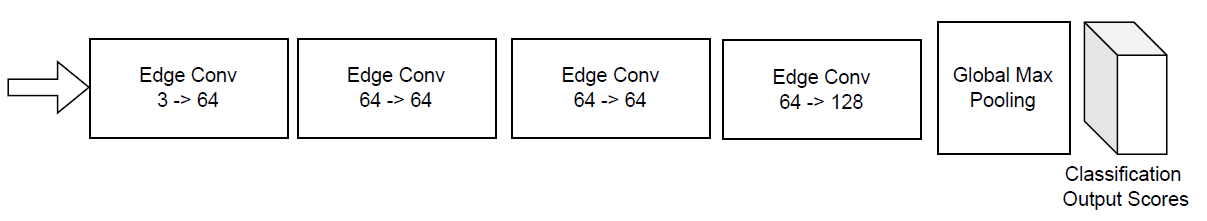}
    \caption{Dynamic Graph Convolutional Neural Network (DGCNN) model architecture.}
		\label{fig:dgcnn}	
\end{figure}

\emph{Dynamic Graph Convolutional Neural Network:}
Instead of working on individual points like PointNet, however, Dynamic Graph Convolutional Neural Network (DGCNN) exploit local geometric structures by constructing a local neighborhood graph and applying convolution-like operations on the edges connecting neighboring pairs of points.
The local neighbourhoods are constructed using the $K$-Nearest Neighbours approach where the top $K$-Nearest Neighbours based on the euclidean distance are considered as local neighbourhood of a point, and edges are drawn from the point to all the points in the local neighbourhood.
The higher the $K$, the more edge dense the graph is.
The KNN connections are computed after every layer.
This is  is why the model is referred to as dynamic.
The DGCNN model architecture is represented in Figure~\ref{fig:dgcnn}.
 
\section{Results and Evaluation}
\label{sec:result}

We consider an acoustic problem within a car compartment.
The objective is to analyze the frequency response function issued from various excitation from the engine.
The results obtained correspond to an acoustic map within the car compartment.
When the noise level becomes too high around the driver, the attention of the driver slow down increasing the risk.
The objective is thus to be able to classify the results into two class: (i) the results where the noise around the driver ear is low; (ii) the results where the noise around the driver is too high.

The acoustic problem is governed by the Helmholtz equation inside the car compartment and boundary conditions defined on the boundary.
Robin boundary conditions modeling absorbing material are defined on the seat of the car and on the roof.
Neumann boundary conditions are defined on the glass of the car.
Dirichlet boundary conditions are defined along the firewall.
The problem is discretized with $P_1$ finite element with stabilized parameter obtained from \cite{HM2004WM}, and the associated linear system of equations is solved with the Alinea library \cite{MC2015IJHPCA}.
The problem is reformulated with a domain decomposition method \cite{MCP2016IJCM} and is solved on Graphics Processing Unit (GPU) \cite{MC2016CC} \cite{CM2017JOS}.
The local solver inside each subdomain is performed with the conjugate gradient method \cite{CM2017AES} involving a auto-tuning of the GPU memory \cite{MCP2015}.
The deep-learning models is implemented with Tensorflow library on GPUs.
The deep-learning models are trained on the dataset divided into train (66.6\%) and test (33.3\%).
A batch size of two objects has been used for the training process for both the models.
The models are trained for 20 epochs each and the accuracy of the model is computed after training on the test set. 

The DGCNN model requires a hyper-parameter, i.e., the $K$-value in $K$-Nearest Neighbours algorithm to be tuned.
To find the most optimal $K $value, experiments are run on different values of $K$ and reported Table~I.
This hyper-parameter tuning step can make a significant difference in final performance of the model.
After considering various values of $K$, the most optimal one is observed to be equal to ten with an accuracy score of 94.5\% and is chosen for the final DGCNN model.

\begin{table}[h!]
\centering
\small
 \begin{tabular}{l r} 
 %\hline
 Value of $K$ parameter & Accuracy\\ [0.5ex] 
 \hline\hline
 5 & 91.67\%\\
 {10} & {94.50\%}\\ 
 15 & 73.61\%\\
 20 & 80.56\%\\
 30 & 80.56\%%\\ [1ex]
 %\hline
 \end{tabular}
\caption{Sensitivity analysis of DGCNN model.}
\end{table}

We can see in Table~II, that the DGCNN model performs well with a classification accuracy of 94.5\%.
The pointnet model, even though it is able to classify the FEM results to certain extent, are out-performed by the DGCNN model.
This might be because of DGCNN's ability to capture local level features much better than the Point-Net model.

\begin{table}[h!]
\centering
\small
 \begin{tabular}{l r} 
 %\hline
 Name of the model & Accuracy\\ [0.5ex] 
 \hline\hline
 CNN(xy projection) & 70.84\%\\
 CNN(yz projection) & 80.00\%\\
 Point-Net model & 79.17\%\\
 DGCNN model & {94.50\%}%\\ 
 %\hline
 \end{tabular}
\caption{Comparison of deep-learning models.}
\end{table}

\section{Conclusions}
\label{sec:conc}

The point-cloud classification models are here considered to classify simulation results arising from finite element analysis, namely the point-net model and the DGCNN model.
Both models work well to automatically analyze the results, and excellent accuracy is obtained using the DGCNN model.

\section{Acknowledgements}
The authors would like to thank M.~Yu, R.~Raj, W.~He, Z.~Liu and Y.~Chen for the useful discussions.


\begin{thebibliography}{1}


\bibitem{Cia2002}
P.G. Ciarlet.
The Finite Element Method for Elliptic Problems.
SIAM.
2002.

\bibitem{Hug2003}
T. Hughes.
The Finite Element Method: Linear Static and Dynamic Finite Element Analysis.
Dover Publications Inc. 2003.

\bibitem{DF2008}
J.E. Dolbow and L.P. Franca.
Residual-free bubbles for embedded Dirichlet problems.
Computer Methods in Applied Mechanics and Engineering. 197(45-48):3751-3759, 2008.

\bibitem{MZ2018CMAME}
F. Magoul\`es and H. Zhang.
Three-dimensional dispersion analysis and stabilised finite element methods for acoustics.
Computer Methods in Applied Mechanics and Engineering, 335:563-583, 2018.

\bibitem{HM2004WM}
I. Harari and F. Magoul\`es.
Numerical investigations of stabilized finite element computations for acoustics.
Wave Motion, 39(4):339-349, 2004.

\bibitem{Rus2006}
A. Russo.
Streamline-upwind Petrov/Galerkin method (SUPG) vs residual-free bubbles (RFB).
Computer Methods in Applied Mechanics and Engineering. 195(13–16):1608-1620, 2006.

\bibitem{HDO1998}
M.A. Hearst, S.T. Dumais, E. Osuna, J. Platt, and B. Scholkopf.
Support vector machines.
IEEE Intelligent Systems and their applications. 13(4):18-28, 1998.

\bibitem{MZE2013EB}
F. Magoul\`es, H.-X. Zhao, and D. Elizondo.
Development of an RDP neural network for building energy consumption fault detection diagnosis.
Energy and Buildings. 62:133-138, 2013.

\bibitem{ZE2012RSER}
H.-X. Zhao and F. Magoul\`es.
A review on the prediction of building energy consumption.
Renewable and Sustainable Energy Reviews. 16(6):3586-3592, 2012.

\bibitem{ZE2012JACT}
H.-X. Zhao and F. Magoul\`es.
Feature selection for predicting building energy consumption based on statistical learning method.
Journal of Algorithms and Computational Technology. 6(1):59-78, 2012.

\bibitem{ZE2010JACT}
H.-X. Zhao and F. Magoul\`es.
Parallel support vector machines applied to the prediction of multiple buildings energy consumption.
Journal of Algorithms and Computational Technology. 4(2):231-250, 2010.

\bibitem{ICN2017}
A. Ioannidou, E. Chatzilari, S. Nikolopoulos, and I. Kompatsiaris.
Deep learning advances in computer vision with 3d data: A survey.
ACM computing surveys (CSUR). 50(2):1-38, 2017.

\bibitem{CM2017JOS}
A.-K. Cheik Ahamed and F. Magoul\`es.
Efficient implementation of Jacobi iterative method for large sparse linear systems on graphic processing units.
The Journal of Supercomputing, 73(8):3411-3432, 2017.

\bibitem{CM2017AES}
A.-K. Cheik Ahamed and F. Magoul\`es.
Conjugate gradient method with graphics processing unit acceleration: CUDA vs OpenCL.
Advances in Engineering Software, 111:32-42, 2017.

\bibitem{MC2016CC}
F. Magoul\`es, A.-K. Cheik Ahamed, and A. Suzuki.
Green computing on graphics processing units.
Concurrency and Computation: Practice and Experience, 28(16):4305-4325, 2016.

\bibitem{MC2015IJHPCA}
F. Magoul\`es and A.-K. Cheik Ahamed.
Alinea: An advanced linear algebra library for massively parallel computations on graphics processing units.
International Journal of High Performance Computing Applications, 29(3):284-310, 2015.

\bibitem{MCP2016IJCM}
F. Magoul\`es, A.-K. Cheik Ahamed, and R. Putanowicz.
Optimized Schwarz method without overlap for the gravitational potential equation on cluster of graphics processing unit.
International Journal of Computer Mathematics, 93(6):955-980, 2016.

\bibitem{MCP2015}
F. Magoul\`es, A.-K. Cheik Ahamed, and R. Putanowicz.
Auto-tuned Krylov methods on cluster of graphics processing unit.
International Journal of Computer Mathematics, 92(6):1222-1250, 2015.

\bibitem{QCH2017}
R.Q. Charles, H. Su, M. Kaichun and L. J. Guibas.
PointNet: Deep Learning on Point Sets for 3D Classification and Segmentation.
In Proceedings of the 2017 IEEE Conference on Computer Vision and Pattern Recognition (CVPR), 77-85, 2017.

\bibitem{SGT2008}
F. Scarselli, M. Gori, A.C. Tsoi, M. Hagenbuchner, and G. Monfardini.
The Graph Neural Network model.
IEEE transactions on neural networks. 20(1):61-80, 2008.

\bibitem{YSL2019}
Y. Wang, Y. Sun, Z. Liu, S.E. Sarma, M.M. Bronstein, and J.M. Solomon.
Dynamic graph CNN for learning on point clouds.
ACM Transactions On Graphics (tog). 38(5):1-12, 2019.

\bibitem{WG2021}
L. Wang and B. Goldluecke.
Sparse-PointNet: See Further in Autonomous Vehicles.
IEEE Robotics and Automation Letters, 6(4):7049-7056, 2021.

\bibitem{SEG2021}
Q. Hu, B. Yang, L. Xie, S. Rosa, Y. Guo, Z. Wang, N. Trigoni, A. Markham.
Learning Semantic Segmentation of Large-Scale Point Clouds with Random Sampling.
IEEE Transactions on Pattern Analysis and Machine Intelligence, 1-1, 2021.

\bibitem{VOX2021}
Y. Li, S. Yang, Y. Zheng and H. Lu.
Improved Point-Voxel Region Convolutional Neural Network: 3D Object Detectors for Autonomous Driving.
IEEE Transactions on Intelligent Transportation Systems, 1-1, 2021.

\end{thebibliography}
\end{document}